\input psfig
\hfuzz=1pt
\font\twelverm=cmr12

\def\Chi{{\bf I}}
\def\ra{\rightarrow}
\def\D{{\cal D}}
\def\H{{\cal H}}
\def\S{{\cal S}}
\def\M{{\cal M}}
\def\Hd{{\cal H}(\Sigma_2)}
\def\Ms{{\cal M}(\Sigma_2, \sigma)}
\def\Rzp{{\cal RO}_0^\prime}
\def\Rz{{\cal RO}_0}
\def\Rp{{\cal RO}^\prime}
\def\R{{\cal RO}}
\def\d{\triangle}
\def\Lua{\Lambda(U,\alpha)}
\def\Lh{\hat{\Lambda}}
\def\Lo{\bar{\Lambda}}
\def\Q{{\bf Q}}
\def\N{{\bf N}}
\def\reals{{\bf R}}
\def\Z{{\bf Z}}
\def\dcup{\sqcup}
\def\ds{d_\ast}
\def\ra{\rightarrow}
\def\QED{  \rlap{$\sqcup$}$\sqcap$ }
\def\Hb{{\bf H}}
\def\Uh{\hat{U}}
\def\Bh{\hat{B}}
\def\cl{\centerline}
\def\tY{\tilde{Y}}
\def\tg{\tilde{g}}
\def\tB{\tilde{B}}
\magnification = \magstep1

\cl{\twelverm Weak Disks of Denjoy Minimal Sets}
\bigskip
\cl{Philip Boyland}
\cl{Institute for Mathematical Sciences}
\cl{SUNY at Stony Brook}
\cl{Stony Brook, NY 11794}
\cl{Internet: boyland@math.sunysb.edu}
\bigskip
{\narrower

 \noindent{\bf Abstract.} This paper investigates the existence of
Denjoy minimal sets and, more generally, strictly ergodic sets
in the dynamics of iterated homeomorphisms. It is shown that 
for the full two-shift, the collection of such 
invariant sets with the weak topology
contains topological balls of all finite dimensions.
One implication is an analogous result that holds
for diffeomorphisms with transverse homoclinic points.
It is also shown that the union of Denjoy minimal sets
is dense in the two-shift and that the set of unique probability
measures supported on these sets is weakly dense in the set of all
shift-invariant, Borel probability measures.

} 
\bigskip

{\bf Section 0: Introduction.} 
One strategy for understanding a dynamical system is to first
 isolate invariant
sets that are  
dynamically  indecomposable. One then studies the structure
of these pieces and how they fit together to give the global dynamics.
This idea goes back 
at least to Birkhoff and has a particularly clear
expression in Conley's Morse decompositions.

There are many notions of dynamical  indecomposibility
in the literature.
In this paper we consider a fairly strong one
that uses both topology and measure.
 An invariant set is called
{\it strictly ergodic} if it is both minimal (every orbit is
dense) and uniquely ergodic (existence of  a unique, invariant
Borel probability measure).
These properties are preserved 
under topological conjugacy but {\it not} measure isomorphism.

The simplest such invariant sets  are periodic orbits, and there
are many theorems concerning their existence.
The next simplest strictly ergodic systems
are probably rigid rotations on the circle with irrational rotation number and
the closely related Denjoy minimal sets. Elements of these
invariant sets  are sometimes 
called (generalized) quasi-periodic points. The models
for  Denjoy minimal sets are the 
minimal sets in  nontransitive circle homeomorphisms with
irrational rotation number. An  abstract dynamical system is
called a {\it Denjoy minimal set} if it is topologically 
conjugate to such a model.
One of the questions that motivated this paper is 
what kind of properties of periodic orbits 
are also true for more general  strictly ergodic invariant  sets,
in particular, for Denjoy minimal sets?

One way to begin to address this question is to 
collect these invariant sets into spaces. For a fixed
homeomorphism $f$ of a compact metric space
$X$, let $\S(X,f)$ denote the set of all strictly ergodic
$f$-invariant subsets of $X$.
Since different minimal sets are of  necessity disjoint,
 each point in $\S(X,f)$ represents
a minimal set that is disjoint from every other  minimal set.
A strictly ergodic set supports a unique invariant Borel
probability measure, so we may use these measures with the
weak topology to put a topology on  $\S(X,f)$.  
 If $\D(X,f)$ denotes the set
of $f$-invariant subsets that are Denjoy minimal sets, then 
$\D(X,f)\subset\S(X,f)$, so we may use the weak topology on
$\D(X,f)$ also.

In ([M]), Mather shows that for a area-preserving monotone
twist map of the annulus, $f:A\rightarrow A$,
 the nonexistence of an invariant
circle with a given irrational  rotation number implies the 
existence of   numerous
Denjoy minimal sets with that rotation number.
 More precisely, using the notation just introduced,
$\D (A, f)$ contains topological balls of every finite dimension.
From one point of view this is a very surprising result. One has 
an arbitrarily large dimensional family of minimal sets
embedded in a two-dimensional dynamical system.
Another question that motivated this paper is how common 
is this kind of phenomenon in dynamics on finite dimensional
manifolds?

It is important to note that even for a smooth system,
$\S(M,f)$ can be empty. One example of this is  Furstenberg's 
 $C^\omega$-diffeomorphism of the two torus
that is minimal but {\it not} strictly ergodic ([F]).
However, for the full shift on two symbols $(\Sigma_2, \sigma)$
one has:

\medskip
{\bf Theorem 0.1.} {\it The space $\S (\Sigma_2, \sigma)$ 
contains a subspace 
homeomorphic to the Hilbert cube and the space
$\D (\Sigma_2, \sigma)$ contains topological balls of
dimension $n$ for all natural numbers $n$. }\medskip

The basic tool in the proof of this theorem is the main  construction.
This construction
takes a certain type of open set in the circle (a regular one) and 
produces a compact, invariant set in the full two-shift.
The construction uses the open set to 
produce itineraries
with respect to a rigid rotation on the circle by an irrational angle. 
This process is somewhat analogous to using a Markov
partition to produce a symbolic model for a system. Another
analogous process is used in the kneading theory of
unimodal maps of the interval. The difference here is that the
chosen open set, in general,  has no relation to the dynamics.
The Hilbert cube of strictly ergodic sets is obtained
by showing that the invariant sets constructed in the two-shift 
have unique invariant probability measures that depend continuously
on the regular open sets in the appropriate topologies.

The main construction is a generalization of Morse and Hedlund's
construction of Sturmian minimal sets as described on page 111
of [G-H]. Such generalizations are a standard tool in 
topological dynamics. In particular, the main construction is a special case 
of the almost automorphic minimal extensions of Markley and 
Paul given in [M-P]. Also of particular relevance are pages 234-241 of [A] 
and [H-H1].

For any regular open set, the main construction
yields a minimal set in the shift. If the open set is a finite
union of intervals, it gives a Denjoy minimal set.
When the open set is more complicated, the resulting minimal set
is more complicated. In particular, it follows from
[M-P] that for certain open sets the construction  gives minimal sets 
that have positive topological entropy and
are not  uniquely ergodic (see Remark 3.4 below).

The full two-shift is frequently embedded in the
iterates of a complicated dynamical system. (In fact, this is one 
definition of a ``complicated'' dynamical system.)
In view of Theorem 0.1 one would therefore expect that
that $\S(X,f)$ will frequently contain a Hilbert cube.
In the following corollary, the first sentence is a consequence
of Theorem 0.1 and the Birkhoff-Smale theorem
 (a particularly suitable statement of which can be found on page
109 of [Rl]).
The second sentence follows from the first and a theorem of
Katok ([K]).
\medskip

{\bf Corollary 0.2.} {\it If $f:M\rightarrow M$ is a diffeomorphism
of the compact manifold $M$ that has a transverse homoclinic orbit
to a hyperbolic periodic point,
 then $\S(M,f)$ contains a subspace 
homeomorphic to the Hilbert cube. In particular,  this is the case 
 when $M$  is two-dimensional, $f$ is $C^{1 + \alpha}$ and has
positive topological entropy. }
\medskip

As was the case with Mather's theorem, one has a large dimensional
family of minimal sets (in this case an infinite dimensional family) embedded
in  finite dimensional dynamics. We shall see in Remark 3.7 below that
in many cases this can be viewed as a manifestation of the 
fact that the Hilbert cube is the continuous, surjective image
of the Cantor set. There is an invariant
Cantor set $\hat{\Lambda}$ embedded in the dynamics.
The orbit closure of each point in  $\hat{\Lambda}$ supports a unique
invariant probability measure.
When the measures are given the
weak topology,
the map that takes the point 
to the  measure 
is a continuous surjection 
of the Cantor set $\hat{\Lambda}$
onto the Hilbert cube.

There are two important examples that illustrate the necessity of
the smoothness and dimension  in the second sentence of Corollary 0.2.
In [R] Rees  constructs
a homeomorphism of the two torus that is minimal
and has positive topological entropy. 
Herman gives a $C^\omega$-difeomorphism
of a $4$-manifold that is also minimal
with positive topological entropy ([Hm]).
 Neither example
is uniquely ergodic, so in these cases $\S(M,f)$ is empty.
The second sentence of Corollary 0.2 also raises the question of
a converse. Specifically, if  $\S(M,f)$ contains a subspace
that is homeomorphic to the Hilbert cube, does $f$ have positive
topological entropy? Proposition 3.1 shows that this is false
on manifolds of  dimension bigger than three.

It is an easy exercise to show that periodic orbits
are dense in the full two-shift. A somewhat deeper result
due to Parthasarathy  says that the invariant probability
measures supported on period orbits are weakly dense in
the set of all shift-invariant probability measures,
$\Ms$ ([P]). The next proposition gives the analog of these
results for Denjoy minimal sets.

\medskip
{\bf Proposition 0.3.}{\it  

{\leftskip=40pt\parindent=-18pt

(a) The set of points that are members of Denjoy minimal sets
is dense in $\Sigma_2$.

(b) The set of invariant measures supported on Denjoy minimal sets
is weakly dense in the set of invariant measures, i.e. $\D(\Sigma_2, \sigma)$
is dense in $\Ms$.

}}
\medskip

This paper is organized as follows. Section 1 gives basic definitions,
background information
and the main construction.
Section 2 contains the statement and proof of the main theorem. This theorem
describes continuity properties of the main construction and  the structure
of resulting invariant sets. Section 2 also contains the proof 
of Theorem 0.3. The proof of Theorem
0.1 is given in Section 3, as is the  example
that shows that the converse of Corollary 0.2  is false in dimensions three
and greater. The last section examines the relationship between the
intrinsic rotation  number of a Denjoy minimal set and its ``extrinsic''
rotation  number when it is embedded in a map of the annulus.
It is also shown that any Denjoy minimal set in the 
two-shift can be generated from a regular open set in the circle
using the main construction. 

\medskip

{\bf Acknowledgments:} The author would like 
to thank B. Kitchens, N. Markley, B. Weiss and S. Williams
for useful comments and references.

\bigskip

{\bf Section 1: Preliminaries.}
This section introduces   assorted notation and  definitions and recalls
some basic  facts from topology, ergodic theory and topological dynamics.
Many of these facts are stated without proof or references. 
In such cases, the facts are either elementary exercises or can be found in
Walters' book [W].

For a set $X, Cl(X), Int(X), X^c$ and $Fr(X)$ denote the
closure, interior, complement and frontier of the set, respectively.
The operator $\dcup$ is the disjoint union. Thus $A\dcup B$ represents the 
union of the two sets, but conveys the added information that the sets
are disjoint. The  indicator function of a set $X$ is denoted
$\Chi_X$. Thus  $\Chi_X(x) = 1$ if $x\in X$, and is $0$ otherwise.
  The circle is $S^1 = \reals/\Z$
 and $R_\eta : S^1\rightarrow S^1$ is rigid
rotation by $\eta$, {\it i.e} $R_\eta(\theta) =
 \theta + \eta\  \hbox{mod} \  1$.
Haar measure on the circle is denoted by $m$.

A nonempty, proper subset $U\subset S^1$  is called a {\it regular open} set
if $Int(Cl(U)) = U$.
The set of all regular open sets is 
$$\R = \{ U\subset S^1 : U \hbox{ is a regular open set} \}.$$ 
Given an open set
$U$, its $\ast${\it -dual} is 
the interior of its complement and is denoted  by $U^\ast = Int(U^c)$.
Note that $U$ is regular open if and only if $S^1$ can be written as
the disjoint union of three nonempty sets, 
$S^1 = U \dcup F \dcup U^\ast$ with $F = Fr(U) = Fr(U^\ast)$.
In consequence, $U\in\R$ if and only if $U^\ast\in\R$.

The set $\R$ of regular open sets will be topologized using
the symmetric difference of sets. For $U, V \in \R$, their symmetric
difference is $U \d V = (U \cap V^c)\dcup (U^c\cap V)$ and 
the distance between them is $d(U,V) = m(U\d V)$.
 If $U$ and $V$ are regular open, when
$U\d V$ is nonempty it contains an interval. In particular,
$d(U,V) = 0$ if and only if $U = V$. Since $d(U,V) = 
\int |\Chi_U - \Chi_V| dm = \|\Chi_U - \Chi_V\|_1$, 
$\R$ maybe thought of as a subspace of $L^1(S^1, m)$. This makes 
it clear that $d$ gives a metric
on $\R$.

If the frontiers of either $U$ or $V$ have 
positive measure, it could happen that $d(U,V) \not = d(U^\ast, V^\ast)$.
To avoid this and related situations it is 
sometimes necessary to restrict attention 
to the set of regular open sets whose frontiers have measure zero,
$$\Rz = \{ U\in \R : m(Fr(U))= 0\}.$$
A metric that controls both regular open sets and their
$\ast$-duals is given by 
$$\ds(U,V) = (d(U,V) + d(U^\ast, V^\ast))/2.$$ 
Unless otherwise
noted, the topology on $\R$ will be that given by the metric $\ds$.
 Note that when restricted
to $\Rz$, $d$ and $\ds$ give the same metric.

It will also be useful to identify regular open sets
that are equal after a rigid rotation of $S^1$. More 
precisely, say $U\sim V$  if there exists an $\eta \in S^1$ with
$V = R_\eta(U)$. Denote the quotient spaces by  $\Rp = \R\;/\!\sim$ and
$\Rzp = \Rz\;/\!\sim$. Note that the topology generated by the projection
$\R\rightarrow\Rp$ can be viewed as being generated by  the metric
$d^\prime([U], [V]) = \inf \{ \ds(U, R_\eta(V)) : \eta\in S^1\}$, where
$[U]$ denotes the equivalence class of $U$ under $\sim$.

A related notion is that of a symmetric set. A set $U\in\R$ is
called {\it symmetric} if there exists an $\eta\not = 0 $ with
$R_\eta(U) = U$. Because $U$ is open, such an $\eta$ will always be
a rational number.

In this paper a {\it dynamical system} means a pair
$(X,h)$ where $X$ is a compact metric space and 
$h$ is a homeomorphism. Given a point $x\in X$, its
orbit is $o(x,h) = \{ \dots, h^{-1}(x), x, h(x), \dots\}$.
A finite piece of the forward orbit is denoted 
$o(x,h, N) = \{ x, h(x), \dots, h^N(x)\}$.
If $(X,h) \ra (Y,g)$ is a continuous semiconjugacy, then 
$(X,h)$ is called an {\it extension} of  $(Y,g)$, and
$(Y,g)$ is a {\it factor} of $(X,h)$. When  the
semiconjugacy is one to one on a dense $G_\delta$
set, the extension is termed {\it almost one to one}.  

The pair $(X,h)$ is called a {\it minimal set} if every
orbit is dense.
The pair is {\it uniquely ergodic} if there
exists a unique invariant Borel probability measure. 
A useful characterization is: $(X,h)$ is uniquely ergodic if and only if
the sequence of functions $(\sum_{i=0}^N f \circ h^i)/ (N +1) $ 
converges uniformly for all $f\in C(X, \reals)$.
A pair that is both minimal and uniquely ergodic is called
{\it strictly ergodic}.  Note that the property of
being minimal, uniquely ergodic
or strictly ergodic is preserved under topological conjugacy.
Also, if an extension is strictly ergodic, then so is its factor.

A compact $h$-invariant set $Y\subset X$ is called  minimal, uniquely ergodic
or strictly ergodic if $h$ restricted to $Y$ has that 
property. In a slight abuse of notation, this situation is described
by saying that $(Y,h)$ is minimal, etc.

Perhaps the simplest nontrivial strictly ergodic system is
$(S^1,R_\alpha)$ for an irrational $\alpha$.
A homeomorphism $g:S^1\rightarrow S^1$ that has an irrational
rotation number and the pair $(S^1, g)$ 
is {\it not} minimal is called a {\it Denjoy
example}. Such examples are classified up to topological
conjugacy in [My]. The two classifying invariants are
the rotation number and the set of orbits that are 
``blown up'' into intervals.  A Denjoy
example always has a unique minimal set $Y\subset S^1$ with
$(Y,g)$ strictly ergodic. 

An abstract dynamical system $(X,h)$ is called a 
{\it Denjoy minimal set} if it is topologically conjugate
to the minimal set in a Denjoy example. Such an $(X,h)$ is
always strictly ergodic. Mather points out
in [M] that a Denjoy minimal set $(X,h)$ always has a well defined
{\it intrinsic rotation number}, {\it i.e.} if $(X,h)$ is topologically
conjugate to the minimal sets in two Denjoy examples $(S^1, g_1)$ and
$(S^1, g_2)$, then either $g_1$ and $g_2$ have the same rotation
number or else $g_1$ and $g_2^{-1}$ do.
If $(X,h)$ is a Denjoy minimal set with intrinsic rotation number
$\alpha$, it is an almost one to one extension of $(S^1,R_\alpha)$.

A general dynamical system $(Z,h)$ can have many invariant subsets
that are Denjoy minimal sets or strictly ergodic. These subsets
are collected together in the spaces
$$\D(Z,h) = \{ Y\subset Z: (Y,h) \hbox{ is a Denjoy minimal
set} \}$$
and
$$\S(Z,h) = \{ Y\subset Z: (Y,h) \hbox{ is strictly ergodic} \}.$$

To topologize these spaces we recall the weak topology
on measures. Given a dynamical system $(Z,h)$, the set of all
its invariant, Borel probability measures is denoted
$\M(Z,h)$. 
The weak topology on $\M$ can be defined by saying that
the measures $\mu_n\rightarrow\mu_0$ weakly if and only if
$\int f d\mu_n\rightarrow\int f \mu_0$ for all 
$f \in C(Z, \reals)$. Note that $\M(Z,h)$ with this
topology is compact, and when viewed as a subspace of the 
dual space to $C(Z, \reals)$, it is convex with extreme points
equal to the ergodic measures. Since a strictly ergodic
system supports a unique invariant probability measure, there is
a natural inclusion $\S(Z,h)\subset\M(Z,h)$. This inclusion 
induces a topology on $\S(Z,h)$ that will be called the weak
topology. The fact that  $\D(Z,h)\subset\S(Z,h)$ allows
us to use the weak topology on $\D(Z,h)$ also.

In the absence of  unique invariant measures  
we use the Hausdorff metric to measure the distance
between compact invariant subsets. Given a compact space $X$, the
space consisting of the closed subsets of $X$ with the Hausdorff topology
is denoted $\H(X)$. Note that if $X$ is compact metric then so is
$\H(X)$. A map $\Phi :E\rightarrow \H(X)$ is called
{\it lower semicontinuous} if for all closed subsets
$Y\subset X$, the set $\{ e\in E : \Phi(e) \subset Y\}$ is closed
in $E$. We will need the fact that the following property implies
that $\Phi$ is lower semicontinuous: When 
$e_n\rightarrow e$ in $E$ and for some subsequence $\{n_i\}$, 
$\Phi(e_{n_i})\ra K$ in $\H(X)$ then $\Phi(e)\subset K$.
Informally, $\Phi$ is lower semicontinuous if when you perturb
$e$, $\Phi(e)$ may get suddenly larger, but never suddenly smaller. 

The {\it full shift on two symbols} is the the pair $(\Sigma_2, \sigma)$
consisting of the sequence space
$ \Sigma_2 = \{0,1\}^\Z$ and the shift map $\sigma$. A {\it symbol block}
$b$ is a finite sequence $b_0, b_1, \dots,  b_{N-1}$ which each
$b_i$ equal to $0$ or $1$. The {\it length} of the block $b$ is $N$
and  the {\it period} 
is its period when considered as a cyclic word.
A sequence $s\in\Sigma_2$ has {\it initial block} $b$ if
$b_i = s_i$ for $i=0, \dots,  N-1$. It is notationally convenient
to view the topology on $\Sigma_2$ as being generated by a metric
$d_\Sigma$ with $d_\Sigma(s, t) < 1/N$ if and only if
 $s_i = t_i$ for $|i|<N$.
A {\it cylinder set} depends on a block $b$ and an integer
$n$ and is a set of the form
$$C^n_b = \{ s\in \Sigma_2 : s_{i + n} = b_i, \hbox{ for } i = 0, 
\dots,  length(b)-1 \}. $$
If $n = 0$, we write $C^0_b = C_b$.

Since cylinder sets are both open and closed, their indicator functions
are continuous. In fact, the finite linear combinations 
 of such indicator functions form a dense set 
in $C(\Sigma_2,\reals)$. This implies that
the measures $\mu_n\rightarrow\mu_0$ weakly if and only if
$\mu_n(C_b^n) \rightarrow\mu_0(C_b^n)$ for all cylinder sets
$C_b^n$. Since the elements of $\Ms$ are shift invariant measures,
any such measure $\mu$ satisfies $\mu(C_b^n) = \mu(C_b)$ for all
$n$. Thus the topology on $\Ms$ is in fact generated by the metric
$$d(\mu_1, \mu_2) = \sum | \mu_1(C_{b^{(n)}}) - \mu_2(C_{b^{(n)}}) | /
2^n$$
where the sum is over some enumeration $b^{(n)}$ of all possible
 blocks by the natural numbers
$n$.

The main construction in this paper takes a regular open
set in the circle and produces a compact invariant set in
$(\Sigma_2, \sigma)$ along with
an invariant measure.  As noted in the introduction, it is closely 
related to the construction given in [M-P]. We are primarily 
interested here in the dependence of the construction on
the open set and a ``rotation number''. This dependence is encoded 
in two functions $\lambda:\Rz \times S^1\rightarrow \Ms$
and $\Lambda:\R \times S^1 \rightarrow \Hd$ defined as follows.

Fix $U\in \R$ and $r\in S^1$  Define $B\subset S^1$ as 
$$B  = \{ x \in S^1 : o(x, R_r) \cap Fr(U) = \emptyset\}.$$
Since $U$ is regular open, $ Fr(U)$ is closed and nowhere dense,
and thus since\break\hfill
 $B = \cap_{i\in \N} R^i_r(Fr(U)^c)$,
$B$ is 
dense $G_\delta$.
 Now define
$\phi: B \ra \Sigma_2$ so that
$$
(\phi(x))_i = \Chi_U(R^i_r(x)).
$$
Thus for any point $x\in B$, the sequence $\phi(x)$ is
the ``itinerary'' of $x$ under $R_r$ with respect to the set
$U$, {\it i.e.} $\phi(x)$ has a $1$ in the i$^{th}$ place
if $R_r^i(x)$ is in $U$ and $0$ if it is in $U^\ast$.
It is easy to see that $\phi$ is continuous. 

Now define
 $\Lambda(U,r) = Cl(\phi(B))$.
If $U\in \Rz$, then $m(U\dcup U^\ast) = 1$ and so $m(B) = 1$. Thus
we may define a probability measure $\lambda\in\Ms$ by
$\lambda = \phi_\ast(m)$, where as usual this means that 
$\lambda(X) = m(\phi^{-1}(X))$ for a Borel set $X$. 

In this construction, $B$ and $\phi$ depend on the choice of
$U$ and $r$. If this dependence needs to be emphasized, we write
$B = B_{U,r}$ and $\phi = \phi_{U,r}$. It is clear that for
all $\eta \in S^1$ and $U\in \R $,
 $\Lambda(R_\eta(U), r) = \Lambda(U,r)$ and
 $\lambda(R_\eta(U), r) = \lambda(U,r)$. 
 Thus the maps $\Lambda$ and
$\lambda$ descend to maps on $\Rp \times S^1$ and $\Rzp \times S^1$
 that will also
be called $\Lambda$ and  
$\lambda$.

To make the last definition, we need to adopt the notation
that $U_0 = U^\ast$ and $U_1 = U$. For a block of symbols
$b$ of length $N+1$, define
$$U_{b,r}=  \bigcap_{i=0}^N R_r^{-i}(U_{b_i}).$$ The important
property of these sets is that for $x\in B$, $x\in U_{b,r}$ if and
only if $\phi(x)$ is in the cylinder set $C_b$. 
As a consequence, for $U\in\Rz$, $\lambda(U,r)[C_b] = \phi_\ast m(C_b) = 
m(U_{b,r})$.

\medskip
{\bf Lemma 1.} {\it The following maps are continuous.

(a) For fixed $U\in\R$, the map $S^1\rightarrow \reals $ given
by $\eta \mapsto \ds(U,R_\eta(U))$.

(b) For fixed $U\in\R$, the map $S^1\rightarrow\R$ given
by $\eta\mapsto R_\eta(U)$.

(c) The map $\R\rightarrow\reals$ given by $U\mapsto m(U)$.

(d) For fixed symbol block $b$, the map
$\R \times S^1 \rightarrow\R$ given by $(U,r)\mapsto U_{b,r}$.
}
\medskip
{\bf Proof of (a) and (b).}
We first prove continuity of the map $\eta \mapsto d(U,R_\eta(U))$ at
$\eta = 0$. Since $U\in S^1$ is open, we can find a countable set of disjoint
intervals $\{I_n\}$ so that $U = \dcup I_n$. Now given $\epsilon > 0$,
pick $M$ so that $\sum_{n>M} m(I_n) < \epsilon/4$ and
assume $|\eta| < \epsilon/ (4 M)$. Now for each $n$, clearly
$m(I_n \cap R_\eta(U)^c) < \eta$ and so
$$\eqalign{m(U \cap R_\eta(U)^c)& < \sum_{n>M} m(I_n) +
 \sum_{n\le M} m(I_n \cap R_\eta(U)^c)\cr
& < \epsilon/2.\cr}$$
Now since $m(U^c \cap R_\eta(U)) = m(R_{-\eta}(U)^c \cap U)$, we
also get $m(U^c \cap R_\eta(U)) < \epsilon/2$ and so $d(U, R_\eta(U))
< \epsilon$.

What we have just shown also 
implies that $\eta \mapsto d(U^\ast,R_\eta(U^\ast))$
is continuous at $\eta = 0$,
and thus $\eta \mapsto d_\ast(U,R_\eta(U))$ is also.
Since $d(R_\eta(U), R_{\eta^\prime}(U)) = d(U, R_{\eta - \eta^\prime}(U))$,
the continuity of $\eta \mapsto R_\eta(U)$ at all $\eta$
follows. Finally, since
$d_\ast$ is a metric, and therefore a continuous function
$\R \times \R \rightarrow \reals$, we get  $\eta \mapsto d_\ast(U,R_\eta(U))$
continuous for all $\eta$.

{\bf Proof of (c).}
Given two finite collections of sets ${A_i}$ and ${B_i}$  
with $i \in \{0, \dots N\}$
using the fact that $d(A, B) = \|\Chi_A-\Chi_B\|_1$ and standard
integral inequalities
it is easy to show that 
$|m(A) - m(B)| \le d(A,B) $ and 
$\ds(\cap A_i, \cap B_i ) \le \sum \ds(A_i, B_i)$. 

The continuity of $U \mapsto m(U)$,
 follows from   the fact that $\ds(U,V) \le \epsilon/2$ implies
 $\epsilon \ge d(U,V) \ge |m(U) - m(V) | $.

{\bf Proof of (d).}
If the length of the fixed block $b$ is $N+1$, then given $\epsilon > 0$
using (a), pick $\delta < \epsilon/(2N + 2)$ so that 
 $|\eta| < \delta$ implies
$d_\ast(U, R_\eta(U)) < \epsilon/(2N + 2)$.

We therefore have for $(V,s)\in \R \times S^1$ with 
$d(U,V) <\delta$ and $|r-s|<\delta/N$,
$$\eqalign{ \ds( U_{b,r}, V_{b,s}) &
    = \ds(\cap R_r^{-i}(U_{b_i}),\cap R_s^{-i}(V_{b_i}))\cr
          &\le \sum d_\ast(R_r^{-i}(U_{b_i}), R_s^{-i}(V_{b_i}))\cr
          &=   \sum\ds(R_{i(s-r)}(U_{b_i}), V_{b_i})\cr
          &\le \sum(\ds(R_{i(s-r)}(U_{b_i}), U_{b_i}) + \ds(U_{b_i},V_{b_i})\cr
           &\le \epsilon.\cr}$$
\line{\hfill \QED}

{\bf Section 2: The main theorem.} The main goal of this section
is to prove the following theorem. For the reader interested in the 
quickest route to Theorem 0.1, we note that the lower 
semicontinuity of $\Lambda$ and the results in part (3) are not needed
for that  proof.
\medskip

{\bf Theorem 2. }{\it Let the maps $\lambda:\Rzp \times S^1\rightarrow \Ms$
and $\Lambda:\Rp \times S^1\rightarrow \Hd$ be as defined in Section 1.

(1) The map $\lambda$ is continuous and the map 
$\Lambda$ is lower  semicontinuous.

(2) Fix $\alpha \not\in \Q$.

{\leftskip=60pt\rightskip=20pt\parindent=-18pt

(a) For all $U\in \Rp$, $(\Lua, \sigma)$ is an almost one to one minimal
extension of 
$(S^1, R_{n \alpha})$ for some natural number  $n$.

(b) If $U\in\Rzp$, then $(\Lua, \sigma)$ is uniquely ergodic.

(c) If $Fr(U)$ is a finite set, then $(\Lua, \sigma)$ is a Denjoy
minimal set with intrinsic rotation number $n \alpha$ 
for some natural number  $n$. 

(d)  For fixed $\alpha \not\in \Q$, when considered as a function of $U$, 
$\Lambda$ and $\lambda$ are injective.

}
(3) Fix $p/q \in Q$ with $p$ and $q$ relatively prime.

{\leftskip=60pt\rightskip=20pt\parindent=-18pt

(a) For all $U\in \Rp$,
$\Lambda(U,p/q)$ is a finite collection of periodic orbits
whose periods divide $q$.

(b) For fixed $p/q \in Q $, when considered as a function of $U$,
the image of $\lambda$ is the convex hull of the probability  measures
supported on the periodic orbits whose periods divide $q$.

}}

\medskip

{\bf Proof of (1).}
Since $\Lambda(R_\eta(U), r) = \Lambda(U,r)$, it suffices to
check the continuity of $\Lambda$ as a map defined on $\R$.
A similar comment holds for $\lambda$.

As noted in the previous section, the weak topology
on $\Ms$ is generated by the metric $d(\lambda_1, \lambda_2)
= \sum | \lambda_1(C_{b^{(n)}}) - \lambda_2(C_{b^{(n)}}) | /
2^n$ and   $\lambda(U,r)[C_b] = m(U_{b, r})$. 
Thus to prove the continuity of $\lambda$ it suffices to check 
that for fixed $b$ the map 
$U \mapsto m(U_{b, r}) $ is continuous. This follows from Lemma 1 (c) and (d).

For the proof of the lower semicontinuity of $\Lambda$, begin
by assuming that $(U^{(n)}, r^{(n)}) \rightarrow (U^{(0)}, r^{(0)})$.
If for some subsequence $\{n_i\}$,
 $\Lambda(U^{(n_i)}, r^{(n_i)}) \rightarrow
K$ in the Hausdorff topology, then we will show that $\Lambda(U^{(0)}, r^{(0)})
\subset K$. As noted in the previous section, this implies the desired 
semicontinuity.
 Fix
an $x_0 \in B^{(0)}$ and integer $N > 0$ and let $b$ 
be the initial block of length $N+1$ in $\phi(x_0)$. This certainly implies
that 
$U^{(0)}_{b, r^{(0)}}$ is a nonempty open set and therefore has positive
 measure.  Therefore by Lemma 1 (c) and (d) there exists an $M$ so that 
$n>M$ implies that $m(U^{(n)}_{b, r^{(n)}}) > 0$. In particular,
for $n>M$, there exists $x_n \in B^{(n)}$ so that
$\phi_n(x_n)$ has its initial block  equal to $b$. There
therefore exits a sequence $x_j \in B^{(j)}$ with 
$\phi_j(x_j) \rightarrow \phi_0(x_0)$.

Now assuming that for some subsequence $\{n_i\}$,
 $\Lambda(U^{(n_i)}, r^{(n_i)}) \rightarrow
K$ in the Hausdorff topology, then if $x_{n_i} \in B^{(n_i)}$ is the
appropriate
subsequence of the sequence constructed in the previous paragraph, 
then  $\phi_{n_i}(x_{n_i}) \rightarrow \phi_0(x_0)$, so certainly
$\phi_0(x_0) \in K$. But $x_0\in B^{(0)}$ was arbitrary, and so
$\phi_0(B^{(0)})\subset K$ and since $\Lambda(U^{(0)}, r^{(0)})$ is 
the closure of the $\phi_0(B^{(0)})$, we have  $\Lambda(U^{(0)}, r^{(0)})
\subset K$, as required.

{\bf Proof of (2).}
For the proof of (2),  fix an $\alpha\not\in \Q$ and for the proof
of (2a), (2b) and (2c) a $U\in\R$. We will suppress the dependence
of various objects on $U$ and $\alpha$ and so $\Lambda = \Lambda(U, \alpha)$,
etc.

{\bf (2a).}  To prove the minimality of 
$\Lambda $ we use the following characterization
of minimality ([O]):
If $f:X \rightarrow X$ is a homeomorphism of a compact 
metric space and $x \in X$, then $Cl(o(x,f))$ is a minimal set if and
only if given $\epsilon > 0$, there exists an $N$ such that for all
$n$, there exists an $i$ with $0 \le i \le N$ and $d(f^{n + i}(x), x) <
\epsilon$.

To apply this to the case at hand, first note that for
$x \in B$, certainly $o(x, R_\alpha)$ is dense in $B$, and so
$\Lambda = Cl(o(\phi(x), \sigma))$. Since $(S^1, R_\alpha)$ is minimal,
the above property holds for $Cl(o(x,R_\alpha))$. Since $\phi$ is 
continuous, it also holds for $Cl(o(\phi(x), \sigma)) = \Lambda$,
which is therefore minimal.

The proof of the semiconjugacy 
requires a new definition.
Given $U, V \in \R$, define $\rho(U,V) = \sup\{ m(I) : I
\hbox{ is an interval contained in\ } U\d V\}$. Now $\rho$ will not satisfy
the triangle inequality but it is easy to see that for fixed $U\in \R$,
the map $\eta \mapsto \rho(U, R_\eta(U))$ is a continuous function
$S^1 \rightarrow \reals$. Also,  if $U$ is asymmetric,
 then $\rho(U, R_\eta(U))=
0$ if and only if $\eta = 0$.

The first step in the proof of the semiconjugacy
is to show that $\phi$ is injective when  $U$ is
asymmetric. 
Assume that for
$x_1, x_2 \in B$, $\phi(x_0) = \phi(x_1)$, and therefore for all $i$,
$\Chi_U(R_\alpha^i(x_1)) = \Chi_U(R_\alpha^i(x_2))$.
 Thus if $x_2 = R_\eta(x_1)$,
$\Chi_U = \Chi_U \circ R_\eta$ when restricted to the dense set
$o(x_1, R_\alpha)$. In particular, $\rho(U, R_\eta(U))=  
0$ and since $U$ is asymmetric, $d(x_1, x_2) = \eta = 0$.

Continuing with the assumption that $U$ is asymmetric, we show
that $\phi^{-1}$ is uniformly continuous. Since $\phi(B)$ is
certainly dense in $\Lambda$, this implies that we can extend $\phi^{-1}$ 
to a semiconjugacy from $(\Lambda, \sigma)$ to $(S^1, R_\alpha)$.

Since $S^1$ is compact and $\eta \mapsto \rho(U, R_\eta(U))$ is 
continuous,  given 
$\epsilon > 0$ there exists a $\delta > 0$ so that $\rho(U, R_\eta(U)) < 
\delta$ implies $|\eta| < \epsilon$. Pick $N >0$ so that
for every $x\in S^1$, every interval of length $\delta$ contains a point of
$o(x, R_\alpha, N)$.
Now if $x_1, x_2 \in B$ satisfy $d_\Sigma(\phi(x_1), \phi(x_2)) < 1/N$
 and if $x_2 = R_\eta(x_1)$, then $\Chi_U = 
\Chi_U \circ R_\eta$ when restricted to the
set $o(x, R_\alpha, N)$. Now if $\rho(U, R_\eta(U)) > \delta$ then
$U\d R_\eta(U)$ will contain an interval of length $\delta$ 
and thus a point of
$o(x, R_\alpha, N)$, a contradiction. Thus $\rho(U, R_\eta(U)) < \delta$
and so by the choice of $\delta$, $d(x_1, x_2) =| \eta | < \epsilon$, proving
the uniform continuity of $\phi^{-1}$. 
Note that $\phi(B)$ is dense $G_\delta$ in $\Lambda$
so the extension is almost one to one.

Now assume that $U$ is symmetric. 
The group of numbers $r$ such that
$R_r(U) = U$ has a rational generator, say $p/q$, with $0 < p/q <1 $
and $p$ and $q$ relatively prime. If $U^\prime = \pi(U)$ where
$\pi:S^1 \rightarrow S^1/R_{p/q}$ is the projection, then 
$\Lambda(U,\alpha)$ has $\Lambda(U^\prime, q \alpha)$ as a 
$q$-fold factor (here we have identified $S^1/R_{p/q}$ with
$S^1$). Since $U^\prime $ is asymmetric,  $\Lambda(U^\prime,q  \alpha)$
has $(S^1, R_{ q\alpha})$ as a factor, finishing the proof
of (2a). 

{\bf (2b).}
Let $\psi$ denote the extension of $\phi^{-1}$ to a continuous
semiconjugacy from $(\Lambda, \sigma)$ to $(S^1,R_{\alpha q})$ and assume that
$m(Fr(U)) = 0$. If $\lambda_1$ and $\lambda_2$ are two
invariant Borel probability measures supported on $\Lambda$, then
since  $(S^1,R_{\alpha q})$ is uniquely ergodic, $\psi_\ast(\lambda_1) = 
\psi_\ast(\lambda_2) = m$. If $X \subset \Lambda$ is a Borel set,
then since $m(B) = 1$, for $i = 1,2$, $\lambda_i(X) = 
\lambda_i(\psi^{-1}(B) \cap X)$. Now since $\psi$ is injective
on $B$, this is equal to $\lambda_i(\psi^{-1}(B \cap\psi(X)) =
m(B\cap\psi(X))$ and so $\lambda_1=\lambda_2$.

{\bf (2c).}
Now
assume $Fr(U)$ is a finite set. In this case, each
$x \in B^c$ will have exactly two preimages under
$\psi$, namely, the limit of $\phi(x_n)$ as $x_n \rightarrow x$
from the right and the limit of $\phi(x_n)$ as $x_n \rightarrow x$ 
from the left. This makes it clear that in this case $\Lambda$ 
is conjugate to the minimal set in the circle homeomorphism obtained
by ``blowing up'' into intervals points on the orbits of each
 $x\in Fr(U)$.

{\bf (2d).}
When $U\in\Rz$, $\Lambda(U,\alpha)$ is the support of $\lambda(U,\alpha)$.
Thus to prove (2d) it suffices to show that $\Lambda(U,\alpha)$ is an injective
function of $U$.
Assume that for some $U_1, U_2 \in \R$,
 $\Lambda(U_1, \alpha) = \Lambda(U_2, \alpha)$.
Using (2a), $\phi(B_1)$ and $\phi_2(B_2)$ are dense $G_\delta$
in the compact metric space $\Lambda(U_1, \alpha) = \Lambda(U_2, \alpha)$.
This implies that $\phi(B_1) \cap \phi_2(B_2)\not = \emptyset$, and so 
there exist $x_1, x_2 \in S^1$ with $\phi_1(x_1) = \phi_2(x_2)$.
Thus if $R_\eta(x_1) = x_2$, then $\Chi_{U_1} = \Chi_{U_2} \circ
R_\eta$ when restricted to the dense set $o(x_1, R_\alpha)$.
This implies that $U_1 \d R_\eta(U_2)$ contains no intervals.
Since the $U_i$ are regular open sets, this means that 
$U_1 = R_\eta(U_2)$ and so $U_1$ and $U_2$ are in the same 
equivalence class in $\Rp$, as required.

{\bf Proof of (3).}
Fix $p/q\in\Q$ with $p$ and $q$ relatively prime.
Since $R^q_{p/q} = Id$, it is clear that any $s \in \Lambda(U, p/q)$
will satisfy $\sigma^q(s) = s$ which implies (3a). Say a 
symbol block $b$ is {\it prime} if its length equals its
period.   For $U\in \Rzp$,
by construction, $\lambda(U, p/q) = \sum m(U_{b, p/q}) \mu_b$ where 
$\mu_b$ is the probability measure supported on the periodic
orbit with repeating block $b$  and the sum is over all prime blocks $b$
whose period divides $q$. With this formula in hand
it is easy to construct a $U$ so that $\lambda(U, p/q)$ is any desired point
in the convex hull given in the statement of (3b). \QED
\medskip
{\bf Proof of Proposition 0.3}
 
(a). A theorem of Parthasarathy  says that the measures supported on
periodic orbits are dense in $\Ms$ ([P]).  Fix one such measure $\mu_0$,
and assume
it is supported on an orbit of period $q$. Using the formula given in
the proof
of Theorem 2 (3b), find a regular open set $U$ with $Fr(U)$ a finite set
and a $p/q$ with $\lambda(U, p/q) = \mu_0$. Now pick irrationals
$\alpha_n\ra p/q$. By Theorem 2 (1), $\lambda(U, \alpha_n)\ra \mu_0$,
and
by Theorem 2 (2c), each $\lambda(U, \alpha_n)$ is the unique measure
supported on a Denjoy minimal set.

(b). It suffices to show that for any symbol block $b$, there
exists an $s\in \Sigma_2$ which has initial block $b$ and
$Cl(o(s, \sigma))$ is a Denjoy minimal set. Fix an irrational
$\alpha$ and $x_0 \in S^1$. Choose a finite
union of intervals $U$ so that $R_\alpha^i(x_0)\in U$ if and only if
$b_i = 1$, for $i = 0, \dots, length(b) - 1$. Further, the open set
$U$ should satisfy $o(x_0, R_\alpha) \cap Fr(U) = \emptyset$. If
$U$ has these properties, Theorem 2 (2c) shows that $\Lambda(U,\alpha)$ is the
desired Denjoy minimal set.\QED

\bigskip

{\bf Section 3: The Hilbert cube of strictly ergodic sets.}
We begin with some definitions in preparation for the
proof of Theorem 0.1.
A copy of the  Hilbert cube  
is given by the collection of sequences,
$$\Hb = \{ \gamma\in\reals^\N : 0\le\gamma_i\le {1\over i+2} \hbox{ for all }
i\in\N\}.$$
A subspace of $\Hb$ that contains topological balls of all
dimensions is
$$\Hb_0 = \{ \gamma\in\Hb : \gamma_i = 0,
 \hbox{ for all but finitely many } i\}.$$
For  $\gamma\in\Hb$, define an asymmetric regular open set $U_\gamma$ by
$$U_\gamma = \bigcup_{i\in\N} ({1\over i+2} - \gamma_i^3,
{1\over i+2} + \gamma_i^3). $$
Now define a map $\Gamma : \Hb \ra \Rzp$ via $\Gamma(\gamma) = [U_\gamma]$.
It is clear that $\Gamma$ is continuous and injective. Since $\Hb$ is 
compact, $\Gamma(\Hb)$ is homeomorphic to $\Hb$.  
\medskip
{\bf Proof of Theorem 0.1.} Fix an irrational $\alpha$.
By Theorem 2 (2ab), the set
$\Lambda(\Gamma(\Hb), \alpha)$ consists of strictly ergodic sets.
Since $\Gamma(\Hb)$ is compact, using Theorem 2 (1) and (2d), we have
that $\lambda(\Gamma(\Hb), \alpha)$ is homeomorphic to $\Gamma(\Hb)$
 and therefore
to $\Hb$. This proves the first statement in the theorem.
To prove the second, note that  Theorem 2 (2c)  implies
that $\lambda(\Gamma(\Hb_0), \alpha)$ consists of measures
supported on Denjoy minimal sets.
 Since $\lambda(\Gamma(\Hb_0), \alpha)$
 is homeomorphic to $\Hb_0$, it (and consequently, $\D(\Sigma_2, \sigma)$)
contains topological balls of all dimensions.\QED
\medskip

{\bf Remarks.}

{\bf (3.1)} In Theorem 0.1 there is an obvious distinction between
$\S(\Sigma_2, \sigma)$, which contains a copy of $\Hb$, and
$\D(\Sigma_2, \sigma)$, which contains a copy of $\Hb_0$.
This is because
$\Lambda(\Gamma(\Hb), \alpha)$
contains minimal sets that are  not Denjoy.
In particular, if $\gamma\in\Hb - \Hb_0$ and for some
$i\neq 0$, $R_\alpha^i(0)\in Fr(U_\gamma)$, then
$\Lambda(U_\gamma, \alpha)$ is not a Denjoy minimal set.
In the semiconjugacy from $(\Lambda(U_\gamma, \alpha), \sigma)$
to $(S^1, R_\alpha)$, the inverse image of $0$ consists of  {\it three}
points.

A Denjoy minimal set is obtained from an irrational
rotation on the circle by replacing (or `blowing up'')
each element of a collection of orbits by a {\it pair} of orbits.
For all $\gamma$ not of the type just described,
 $\Lambda(U_\gamma, \alpha)$ is a Denjoy minimal set.
When $\gamma \in \Hb_0$, the number of orbits blown up
is the same as the number of distinct orbits containing points
of $Fr(U_\gamma)$. For $\gamma\in\Hb - \Hb_0$, if for all
$i\neq 0$, $R_\alpha^i(0)\not\in Fr(U_\gamma)$, 
then $\Lambda(U_\gamma, \alpha)$ is a Denjoy minimal set with
countably many orbits blown up.  All the infinite
dimensional families we  could construct had the property that
some minimal set was not Denjoy. 

{\bf (3.2)} Morse and Hedlund's construction of Sturmian minimal
sets corresponds to the special case $U = (0, \alpha)$.
In this case, $\Lambda(U, \alpha)$ is a Denjoy minimal
set with a single orbit blown up.

{\bf (3.3)} Theorem 2 (1) states that $\Gamma$ is a lower semicontinuous
function whose range is the set of closed subsets
of a compact metric space. When such functions have a domain
that is a Baire space, they are continuous on a dense, $G_\delta$
set (see page 114 of [C]). 
It seems unlikely that $\R$ is a Baire space,
but since $\Gamma(\Hb)$ is
homeomorphic to the Hilbert cube, we may apply this result to
show that the map (for fixed $\alpha$)
$$\Lambda(\cdot\  , \alpha) : \Gamma(\Hb)\ra\Hd$$	
is continuous at a generic point of $\Lambda(\Hb)$.
This result can also be obtained directly by showing that
the map is, in fact, continuous at all points 
$\Gamma(\gamma)$ for which all points of 
$Fr(U_\gamma)$ are on disjoint orbits.

{\bf (3.4)} As is perhaps obvious from Remark (3.1),
when $Fr(U)$ is more complicated topologically, 
so is the structure of $\Lambda(U, \alpha)$ (for 
irrational
$\alpha$). However, Theorem 2 (2b) says that for
all $U\in \R_0$, $\Lambda(U, \alpha)$
is uniquely ergodic.  It is in fact  measure 
isomorphic to
$(S^1, R_\alpha)$. To get minimal sets with more 
interesting
measure theoretic properties we must have 
$m(Fr(U)) > 0$.
In this case the set $B_{U,\alpha}$ from the main 
construction
is a zero measure, dense $G_\delta$ set in the circle.
This leads one to  expect that  $\Lambda(U, \alpha)$ 
could support
more than one invariant probability measure.  

The results of  [M-P] show that this is frequently the 
case.  The relevant construction from that paper 
begins with a Cantor  $K$ in the circle. 
The complement of $K$ is the disjoint union of open 
intervals. One chooses a set of labels for these open 
sets with each open set labeled by zero or one. The  
set of labels is used to construct 
a minimal set in the two-shift as in the main 
construction. If $K$ has positive measure, then for 
most  sets of labels (in the appropriate sense) the 
constructed minimal set is not uniquely ergodic and 
has positive topological entropy.

However, the constructed minimal set can be uniquely 
ergodic as the following example suggested by  Benjamin Weiss 
shows. 
 Let $(X, f)$ be a Denjoy minimal
 set with intrinsic rotation number $\alpha$. Note that $(X, f)$
is both measure isomorphic to  and 
an almost one to one extension of $(S^1, R_\alpha)$. Using  
results of Jewett and Kreiger 
we may find
a zero-dimensional strictly ergodic system $(Z,h)$
that is mixing and has positive  entropy.  
  Let $(Y, g)$ be the product of the two systems. 
Because $(Z,h)$ and  $(X,f)$  are
strictly ergodic and $(Z,h)$ is  mixing and $(X,f)$ has
pure point spectrum,  $(Y,g)$
is strictly ergodic.

Now  think of Y as an extension of X. The main theorem 
and the  remark 
following Theorem 4 in [F-W] imply that 
there is a minimal almost 1-1 extension of X , say 
$(\tilde{Y}, \tilde{g})$,   which maps onto $(X,f)$ 
in such a way   that the invariant 
measures of $(\tilde{Y}, \tilde{g})$ are in one to 
one correspondence with the $g$-invariant measures on  
$Y$. Thus  $(\tilde{Y}, \tilde{g})$ is a strictly 
ergodic, positive entropy,  almost 1-1 
extension of rotation by alpha. Further, as a 
consequence of the method of construction in [F-W],
since $X$, $Y$, and $Z$ are  zero-dimensional, $\tilde{Y}$ is   
also.

 Let $p:\tY\ra S^1$ denote the given semiconjugacy and 
 let $\tB \subset\tY$ be the dense $G_\delta$ set on 
 which $p$ is injective. Pick two sets, each open and 
 closed, with $V_0 \dcup V_1 =\tY$. Note that $U = 
 (p(V_0))^c$ is a regular open set. Use the  partition 
 $\{ V_0, V_1\}$
 in the usual way to get a symbolic model by defining
 $k:\tY\ra\Sigma_2$ so that
 $$(k(y))_i = \Chi_{V_1}(\tg^i (y)).$$
 It is fairly straightforward to show that $\tB\subset 
 p^{-1}(B_{U,\alpha})$ and thus,  $p = \psi\circ k$ where 
 $\psi:\Lambda(U,\alpha)\ra S^1$ is the semiconjugacy 
 constructed in the proof of Theorem 2 (2a). This 
 implies that $\Lambda(U,\alpha)$ is a factor of
 $(\tY,\tg)$, and thus is strictly ergodic. Further, we may choose
 $V_0$ and $V_1$ so that $\Lambda(U,\alpha)$ has
 positive entropy. To finish, note 
that  $m(Fr(U)) > 0$, for if  
 not, $\Lambda(U,\alpha)$ would be measure isomorphic 
 to the zero entropy system $(S^1, R_\alpha)$.

It would be interesting to have conditions on a
regular open set with positive measure frontier that 
distinguish these two cases. More precisely, give 
necessary and sufficient conditions for the unique 
ergodicity of  $\Lambda(U, R_\alpha)$.
 Another interesting question
is the structure
of the set of its invariant measures 
in the cases when $\Lambda(U,\alpha)$ is not uniquely
ergodic ({\it cf.} [Wm]).

{\bf (3.5)} Since each point in  $\S(X,f)$ represents a disjoint minimal
set, the size of $\S(X,f)$ should give some indication of the 
complexity of the dynamics of $f$. The topological entropy of
$(X,f)$, denoted $h(X,f)$,  
is perhaps the most common way of measuring dynamical complexity.
Corollary 0.2 shows that, at least in some cases, when the
topological entropy is positive, $\S(X,f)$ is large. 
If the size of $\S(X,f)$ is to give a measure of dynamical complexity,
the converse should be true. The next proposition shows that this is
not the case, at least when the ``size'' of $\S(X,f)$  is measured by
the maximal dimension of an embedded ball and $X$ is a manifold
of dimension greater than two.

\medskip
{\bf Proposition 3.1.} {\it

 \leftskip=40pt\parindent=-18pt

(a) There exists a compact shift invariant set $\Lh\subset\Sigma_2$
such that $\S(\Lh, \sigma)$ is homeomorphic to the Hilbert cube
and $h(\Lh, \sigma) = 0$.

(b) On any smooth
manifold $M$ with dimension greater than two there exists
a $C^\infty$ diffeomorphism $f$ such that  $h(f) = 0$ and $\S(M,f)$
contains a subspace homeomorphic to the Hilbert cube.

}
\medskip
{\bf Proof of (a).}  Fix an irrational $\alpha$ and let
$T = S^1 \times \Hb$. Define $F:T\ra T$ as $F = R_\alpha\times Id$. 
We will do a construction  analogous to the main construction, but now
using the space $T$ and the map $F$. To get  an open set in $T$
we use  the open sets $U_\gamma$ constructed above
to define $$\Uh= \bigcup_{\gamma\in \Hb} U_\gamma \times\{ \gamma\}. $$
Next let 
$$\Bh = \{ \beta\in T : o(\beta, F) \cap Fr(\Uh) = \emptyset \}$$
and define $\Phi: \Bh\ra \Sigma_2$ so that
$$(\Phi(\beta))_i = \Chi_{\Uh}(F^i(\beta)).$$
Finally, let $\Lh = Cl(\Phi(\Bh))$.

Note that for fixed $\gamma$, $\Phi$ restricted to 
$(S^1\times\{ \gamma\})\cap\Bh$ is just $\phi_{U_\gamma, \alpha}$ from the
main construction and that
$$\Lh = Cl(\bigcup_{\gamma\in \Hb} \Lambda(U_\gamma, \alpha)). $$

Theorem 2 (2a) and (2d) imply that $\Phi$ is injective. Using an argument
similar to one in the proof of Theorem 2 (2a), one gets that $\Phi^{-1}$ is
uniformly continuous, and therefore has a continuous extension to a 
$\Psi:\Lh\ra T$ that satisfies $\Psi\circ\sigma = F\circ\Psi$.

The variational principle
(see page 190 in [W]) implies that $h(\Lh, \sigma) = 0$
if all ergodic measures for $(\Lh, \sigma)$ have metric entropy zero.
If $\eta$ is an ergodic, invariant Borel  probability measure for 
$\Lh$, then $\Psi_\ast(\eta)$ is such a measure for $(T, F)$ and so
$\Psi_\ast(\eta)$ is Haar measure on $S^1\times\{\gamma_0\}$ for some
$\gamma_0$. This implies that $\eta$ is 
supported on $\Psi^{-1}(S^1\times\{\gamma_0\})$.
 Once again, using an argument virtually identical to one in the proof of 
Theorem 2 (2b), one obtains $\eta = \lambda(U_{\gamma_0}, \alpha)$. This
measure with the shift is measure isomorphic to rotation on the circle
by $\alpha$ and therefore has zero metric entropy, as required.
Note that the argument just given also shows that $\S(\Lh, \sigma)$ is
in fact homeomorphic to the Hilbert cube, $\Hb$.

{\bf Proof of (b).}
We first construct the map on the space $P=D^2\times [-1,1]$, where
$D^2$ is a closed two-dimensional  disk.
  Let $h:D^2\ra D^2$ be a Smale horseshoe, {\it
i.e.} $h$ is a $C^\infty$-diffeomorphism whose nonwandering set
consists of the union of a finite number of fixed points
and a set $\Omega$  on which the dynamics are conjugate to
the full two-shift. The compact invariant set $\Lh$ constructed in
the proof of (a) is embedded in $\Omega$ by the conjugacy. Call
this embedded set $\Lo$.

Next, let $h_t$ for $t\in [-1,1]$ be an isotopy with $h_{-1} = Id$,
$h_0=h$, and $h_1 = Id$. Further, $h_t$ restricted to the boundary of $D^2$
should be the identity for all $t$.
 Now pick a $C^\infty$-function $w:P\ra\reals$ with $w\ge 0$ and
$w^{-1}(0) = \partial P \dcup (\Lo \times \{0\})$. Let $g:P\ra P$ be the
time one map of the flow generated by the vector field 
$w(u) {\partial\ \over\partial z}$, where $u = (x,y,z)$ is a point in $P$.
Now let $f = g  \circ (h_t\times Id)$.
 By construction, the nonwandering set
of $f$ is $\partial P \dcup (\Lo \times \{0\})$ and thus $h(f) = 0$. Since
each point on $\partial P$ is a fixed point for $f$, $\S(P,f)$ is homeomorphic to
$\S(\Lo,\sigma)\dcup\partial P$, which in turn, is homeomorphic to
$\Hb \dcup \partial P$.

To obtain the result on a general manifold of dimension three
or higher, embed a copy of $(P,f)$ in it and extend $f$ by the identity
on the rest of the manifold.\QED
\medskip

{\bf Remarks}

{\bf (3.6)} This proposition leaves open the possibility of a converse
to Corollary 0.2 in dimension 2. In this dimension there are a number
of results that show that the existence of certain types of
zero entropy invariant sets can imply that a homeomorphism has
positive topological entropy. For example, if an orientation-reversing
homeomorphism of a compact surface of genus $g$ has periodic orbits with
$g+2$ distinct odd periods, then it has positive entropy ([B-F], [H]).
 For orientation-preserving homeomorphisms there
are restrictions on the periods that occur in zero entropy maps given in 
[S].  Even a single period orbit can imply positive entropy
if the isotopy class on its complement is nontrivial ([Bd]).
These results give credence to the conjecture that for a manifold
$M$ of  dimension $2$, if $f:M\ra M$ is a homeomorphism and 
$\S(M,f)$ contains a topological ball of dimension $3$, then
$h(f)>0$.

{\bf (3.7)} It was noted in the introduction that the existence of a Hilbert
cube of strictly ergodic sets can often be viewed as a manifestation of
a standard topological fact, namely, the Hilbert cube is the
continuous surjective image of the Cantor set. For concreteness,
let $f: M\ra M$ be a homeomorphism with an invariant set $\Lo$ with
$(\Lo, f)$ conjugate to $(\Lh, \sigma)$, where $\Lh$ is the set
constructed in the proof of Proposition 3.1 (b). Using the conjugacy,
the proof of Proposition 3.1 (b), and Theorem 2 (2b)
 one gets that for each $x\in\Lo$,
$Cl(o(x,f))$ supports a single invariant probability measure which is  
$c_\ast(\lambda(U_{\gamma(x)}, \alpha))$ for the appropriate
$\gamma(x)$. Further, the map $x\mapsto c_\ast(\lambda(U_{\gamma(x)}, \alpha))$
is continuous. (More formally, this map is
$$x\mapsto c_\ast(\lambda(\Gamma(\pi_2(\Psi(x))), \alpha))$$
where $\pi_2: S^1\times\Hb\ra\Hb$ is the projection).
The domain of this map is the invariant Cantor set $\Lo$ and
its image is $\lambda(\Gamma(\Hb),\alpha)$, which is homeomorphic
to the Hilbert cube, $\Hb$.

{\bf (3.8)} The construction in the proof of 3.1 (b) can be used
to embed any compact shift invariant subset of $\Sigma_2$ as
the only ``interesting'' dynamics in a three-dimensional diffeomorphism. 
It  is reminiscent
of Schweitzer's construction of 
C$^1$-counterexample to the Seifert conjecture ([Sc]).

\bigskip

{\bf Section 4: Intrinsic and extrinsic rotation numbers.}
In the Section 1 it was noted that abstract Denjoy minimal
sets have  well-defined intrinsic rotation numbers.
The next proposition specializes  some previous
results to the case of fixed intrinsic rotation number.

{\bf Proposition 4.1.} {\it Fix an irrational $\alpha$ 
and let $\D_\alpha(\Sigma_2, \sigma)$ denote the set
of Denjoy minimal sets in the shift with intrinsic rotation
number $\alpha$. 
 
{\leftskip=40pt\parindent=-18pt

(a)  When given the weak topology,
the space $\D_\alpha(\Sigma_2, \sigma)$ contains
topological balls of
dimension $n$ for all natural numbers $n$.
 
(b) The set of points that are members of Denjoy minimal sets
with intrinsic rotation number $\alpha$ is dense in $\Sigma_2$.
 
(c)   If $(D,\sigma)$ is a Denjoy minimal set
 with intrinsic rotation number $\alpha$, then
$ D = \Lambda(U,\alpha)$ for some regular open
set $U$ with $m(Fr(U)) = 0$. Consequently,
$\D_\alpha(\Sigma_2, \sigma)$ $ \subset \lambda(\Rz, \alpha)$.
 
}}

\medskip
{\bf Proof of  Proposition 4.1.}
When  $U$ is asymmetric and $\Lambda(U, \alpha)$ is a
Denjoy minimal set, it has  intrinsic rotation number $\alpha$.
This follows from Theorem 2 (2b) (and its proof).
Thus to prove (a) we need only note that  the proof of Theorem 0.1
 began  with a statement, `` Fix an
irrational $\alpha$''. 
The proof of 
Proposition 0.3 (b) contains a similar statement, so that proof 
proves (b).

To prove (c), note that by definition,
 there exists a conjugacy $c:D\ra Y$ where $Y$ is the minimal set
in a Denjoy example $g:S^1\ra S^1$ with rotation number $\alpha$.
It is a standard fact that there
exists a semiconjugacy $h$ of $(S^1, g)$ to $(S^1, R_\alpha)$
with the properties that $h$ is injective on a set that is dense in $Y$ and
the lift of $h$ is weakly order preserving, {\it i.e.} 
$x < y$ implies $\tilde{h}(x)\le \tilde{h}(y)$.

Now let $p = h\circ c$ and 
$U=(p(C_0))^c$. Since $C_0$ is compact in $\Sigma_2$, $U$
is open. Further,
the properties given above imply that 
$U^\ast = (p(C_1))^c$ and $p(C_0)\cap p(C_1) = Fr(U) = Fr(U^\ast)$.
Thus using a fact from Section 1, $U$ is a regular open set, and 
by construction, $\Lambda(U,\alpha) = D$. Since $p(C_0)\cap p(C_1)$ is
at most countable, $m(Fr(U)) = 0$. \QED

These results, of course, also hold for homeomorphisms with a full
two-shift embedded in their dynamics. In this case, however, one is
perhaps more interested in {\it extrinsic} properties of invariant
sets, {\it i.e.} properties associated with how the sets are 
embedded in the manifold. Perhaps the simplest such extrinsic 
property is the extrinsic rotation number, and the simplest case
in which this can be defined is for a homeomorphism of
the annulus.

If $f:A \ra A$ is a homeomorphism of the annulus and $z\in A$, define
the rotation number of $z$ under $f$ as 
$$\rho(z) = \lim_{n\ra\infty} 
{\pi_1(\tilde{f}^n(\tilde{z})) - \pi_1(\tilde{z}) \over n},$$
if the limit exists. Here $\tilde{f}:\reals\times [-1,1]\ra\reals\times [-1,1]$
and $\tilde{z}$ are lifts of $f$ and $z$, respectively,  and 
$\pi_1:\reals\times [-1,1]\ra\reals$ is the projection. 
Note that the rotation number is only defined modulo $1$ as it
depends on the choice of lift.

If $D\subset A$ is a Denjoy minimal set under $f$, then
it is uniquely ergodic. Thus  for all $z\in D$,
$\rho(z)=\int r(z)\; d\mu$, where
$\mu$ is the the unique invariant probability  measure of
$(D,f)$ and $r:S^1\ra \reals$ is the map that lifts 
to $\pi_1\circ\tilde{f} - \pi_1$.
This number will be called the {\it extrinsic rotation
number} of $(D,f)$.

The Denjoy minimal sets constructed by Mather in [M] have
monotonicity properties that imply that their
extrinsic and intrinsic rotation numbers are rationally related.
For Denjoy minimal sets in a general homeomorphism of the annulus 
this will not be the case. As a specific example, we will consider
homeomorphisms $f: A\ra A$ that have a {\it rotary horseshoe}
({\it cf.} [H-H2])
 A picture of the lift of such a map is shown
in Figure 1. The dotted vertical lines are the boundaries of fundamental
domains.
\midinsert
\centerline{\psfig{figure=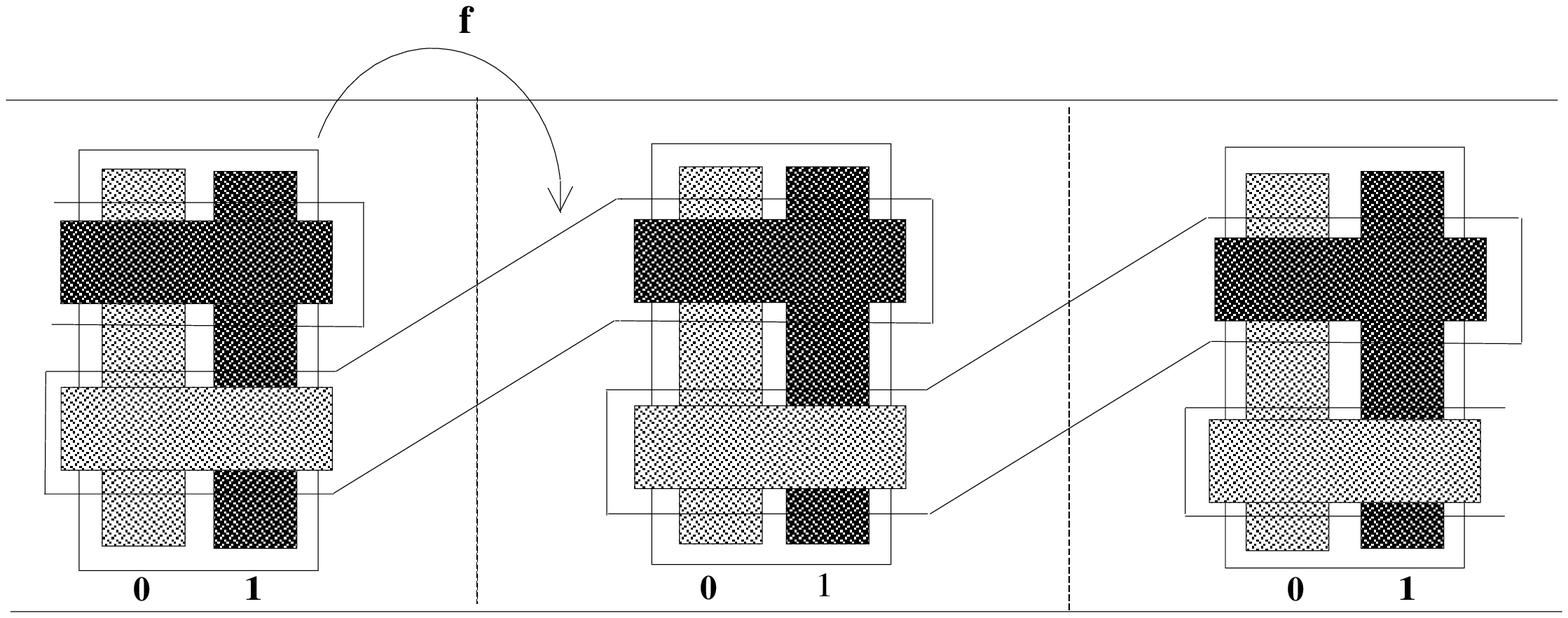,height=2in}} 
\medskip
\centerline{Figure 1: The lift of a rotary horseshoe.}
\endinsert

 A map contains  a rotary horseshoe if it has a 
compact invariant set $\Omega$ that is conjugate to the full two-shift.
The conjugacy $c:\Omega\ra \Sigma_2$ is required to have the property that
for $z\in \Omega$  the first element in
$c(z)$ is $1$ if and only if $\tilde{f}$ moves $\tilde{z}$
(approximately) one fundamental domain to the right. More precisely,
for $z\in\Omega$ it is required that
$$\rho(z) = \lim_{N\ra\infty}
\sum_{i=0}^N {\Chi_{C_1}(\sigma^i(c(z)))\over (N+1)}.$$
Thus $\rho(z)$ is the asymptotic average number of ones 
in the sequence $c(z)$.

We are now almost in a position to state a result about the existence of
Denjoy minimal sets with given intrinsic and extrinsic rotation
number.  For an annulus homeomorphism $f$, let
$\D_{\alpha, \beta}(A,f)$ denote the set of  all  Denjoy
minimal sets for $f$ with intrinsic rotation number $\alpha$
and extrinsic rotation number 
 $ \beta$.
\medskip
{\bf Proposition 4.2.} {\it If a homeomorphism $f:A\ra A$ has a rotary
horseshoe, then for all irrational $\alpha$, and all $\beta\in S^1$, 
$\D_{\alpha, \beta}(A,f)$ contains topological balls of
dimension $n$ for all natural numbers $n$.}
\medskip

{\bf Proof of Proposition 4.2.} If for a given $U\in\Rz$ and irrational
$\alpha$, $\Lambda(U,\alpha)$ is a Denjoy minimal set, then the comments above
Lemma 1 and unique ergodicity imply that for all $s\in \Lambda(U,\alpha)$,
$$\lim_{N\ra\infty}\sum_{i=0}^N {\Chi_{C_1}(\sigma^i(s))\over
 (N+1)}= \lambda(U,\alpha)[C_1]= m(U).$$
 This implies that the corresponding Denjoy minimal set in the
annulus has extrinsic rotation number equal to $m(U)$.
To finish the proof, one need only imitate the proof of Theorem 0.1
using a family $U_\gamma$ that satisfies $m(U_\gamma) = \beta$, for
all $\gamma$.\QED
\medskip

 Note that the case of rational $\beta$ is included in this
result. This means that large dimensional balls of Denjoy minimal sets
with a given {\it rational}
 extrinsic rotation number are present in the dynamics.

\bigskip

\centerline{REFERENCES}
\medskip

\item{[A]} Auslander, J., {\it Minimal Flows and their Extensions},
North Holland Mathematics Studies, vol. 153, 1988.

\item{[B-F]} Blanchard, P. and Franks, J., The dynamical complexity
of orientation reversing    
homeomorphisms of surfaces, {\it Inv. Math.}, {\bf 62}, 1980,  333--339.

\item{[Bd2]} Boyland, P. , An analog of Sharkovski's theorem for twist maps,
 {\it Contemporary Math.}, {\bf 81}, 1988, 119--133.

\item{[C]} Choquet, G., {\it Lectures on Analysis}, vol 2,
W.A. Benjamin, 1969.

\item{[F]} Furstenberg, H., Strict ergodicity and transformations
of the torus, {\it Amer. J. Math.}, {\bf 83}, 1961,  573--601.

\item{[F-W]} Furstenberg, H. and Weiss, B., On almost 
1--1 extensions, {\it Isr. J. Math.}, {\bf 65}, 1989,
311--322.

\item{[GH]} Gottschalk, W. and Hedlund, {\it Topological Dynamics}, AMS Colloquium
Pub., vol. 36, 1955.

\item{[H]} Handel, M. , The entropy of orientation reversing 
homeomorphisms of surfaces, {\it Topology}, {\bf 21}, 1982,  291--296.

\item{[Hm]} Herman, M., Construction d'un diff\'eomorphisme
d'entropie topologique non nulle,
 {\it Ergod. Th. \& Dynam. Sys.},
{\bf  1}, 1981, 65--76.

\item{[H-H1]} Hockett, K. and Holmes. P., Bifurcation to rotating
Cantor sets in maps of the circle, {\it Nonlinearity}, {\bf 1},
1988,  603--616.

\item{[H-H2]} Hockett, K. and Holmes, P., Josephson's 
junction, annulus maps, Birhoff attractors, horseshoes 
and rotation sets, {\it Ergod. Th. \& Dynam. 
Sys.}, {\bf 6}, 1986, 205--239.

\item{[K]} Katok, A., Lyapunov exponents, entropy
and periodic orbits for diffeomorphisms, {\it Publ Math IHES}, {\bf  51},
 1980, 137--173.

\item{[M]} Mather, J., More Denjoy minimal sets for area-preserving
mappings, {\it Comm. Math. Helv.}, {\bf 60}, 1985, 508--557.

\item{[My]} Markley, N., Homeomorphisms of the circle without 
periodic orbits, {\it Proc. Lond. Math. Soc.}, {\bf 20},  1970,
688--698.

\item{[M-P]} Markley, N. and Paul, M., Almost 
automorphic symbolic minimal sets without unique 
ergodicity, {\it Isr. J. Math.}, {\bf 34}, 259--272.

\item{[O]} Oxtoby, J., Ergodic sets,  {\it Bull. A.M.S.}, {\bf 58}, 1952,
116--136.

\item{[P]} Parthasarathy, K.,  On the category of ergodic
measures, {\it Ill. J. Math}, {\bf 5}, 1961.

\item{[R]} Rees, M., A minimal positive entropy
homeomorphism of the two torus, {\it J. London Math Soc}, {\bf 23}, 1981,
537--550.

\item{[Rl]} Ruelle, D., {\it Elements of Bifurcation Theory}, Academic
Press, 1989.

\item{[Sc]} Schweitzer, P. A., Counterexamples to the Seifert conjecture and
opening closed leaves of foliations, {\it Ann. of Math.}, {\bf 100},
1974, 386--400.

\item{[S]} Smillie, J.,  Periodic points of surface homeomorphisms
with zero entropy,
{\it Ergod. Th. \& Dynam. Sys.},
{\bf 3}, 1983, 315--334.

\item{[W]} Walters, P., {\it An Introduction to Ergodic Theory}, Graduate
Texts in Mathematics, vol. 79,  Springer-Verlag,
1982.

\item{[Wm]} Williams, S., Toeplitz minimal flows that are not 
uniquely ergodic,
{\it Z. Wahrscheinlichkeitstheori verw. Geb.}, {\bf 67}, 1984, 95--107.

\bigskip
\line{\hfill revised 2/3/93}

\bye